\documentclass[11pt]{article}
\usepackage{enumerate}
\usepackage{epsfig,graphicx,graphics,latexsym,amssymb,amsfonts,amsmath, amscd}
\usepackage{changebar}
\usepackage{Vmargin}
\usepackage{graphicx}
\usepackage{makeidx}
\usepackage{amsmath}
\usepackage{amssymb}
\usepackage{euscript}
\newtheorem{definition}{{\bf Definition}}[section]
\newtheorem{theorem}[definition]{{\bf Theorem}}

\newtheorem{corollary}[definition]{{\bf Corollary}}
\newtheorem{proposition}[definition]{\noindent {\bf Proposition}}
\newtheorem{lemma}[definition]{\noindent {\bf Lemma}}
\newtheorem{fact}[definition]{\noindent {\bf Fact}}

\newtheorem{claim}[definition]{\noindent {\bf Claim}}

\newtheorem{questions}[definition]{\noindent {\bf Questions}}

\newtheorem{problem}[definition]{\noindent {\bf Problem}}

\newtheorem{remark}[definition]{\noindent {\bf Remark}}

\def\Proof{{\parindent0pt {\bf Proof.\ }}}

\def\endproof{\hfill {\kern 6pt\penalty 500
   \raise -0pt\hbox{\vrule \vbox to5pt {\hrule width 5pt
  \vfill\hrule}\vrule}}}
\newcommand{\K}{\mathbb{K}}
\newcommand{\N}{\mathbb{N}}

\newcommand{\Q}{\mathbb{Q}}

\newcommand{\s}{\scriptstyle}

\begin{document}

\title{Hypomorphy of  graphs up to complementation}
\author{Jamel Dammak \\
\\
 { \it D\'epartement de Math\'ematiques, Facult\'e des Sciences de Sfax}\\
{ \it B.P. 802, 3018 Sfax, Tunisie}\\
E-mail : jdammak@yahoo.fr \\
\\
G\'erard Lopez\\
\\
 { \it Institut de Math\'ematiques de Luminy,
CNRS-UPR 9016,}\\
{ \it 163 avenue de Luminy, case 907, 13288 Marseille
cedex 9, France}\\
E-mail : gerard.lopez1@free.fr  \\
\\
Maurice Pouzet and   Hamza  Si Kaddour\\
\\
 { \it PCS, Universit\'e Claude Bernard Lyon1,}\\
{ \it Domaine de Gerland - b\^at. Recherche B }\\
{ \it  50 avenue Tony-Garnier}\\
{\it  F 69366 - Lyon cedex 07, France}\\
E-mail : Maurice.Pouzet@univ-lyon1.fr ,  Hamza.Sikaddour@univ-lyon1.fr }
\maketitle
\thanks{This  work  was partially supported by CMCU.}

\begin{abstract}Let $V$ be a set of cardinality $v$ (possibly infinite). Two graphs $G$ and $G'$ with vertex set $V$ are  {\it isomorphic up to complementation} if $G'$ is isomorphic to $G$ or to the complement 
$\overline G$ of $G$. Let $k$ be a non-negative integer, $G$ and $G'$ are {\it $k$-hypomorphic  up to complementation} if for every $k$-element subset $K$ of $V$, the induced subgraphs $G_{\restriction K}$ and $G'_{\restriction K}$  are isomorphic up to complementation.  A graph $G$ is {\it $k$-reconstructible up to complementation} if every  graph  $G'$ which is 
 $k$-hypomorphic to   $G$ up to complementation is in fact isomorphic to $G$ up to complementation. 
  We give a partial characterisation of the set
$\mathcal S$ of pairs $(n,k)$ such that  two graphs $G$ and $G'$ on the same set of $n$ vertices are equal up to complementation whenever they are $k$-hypomorphic up to complementation. We prove in particular that $\mathcal S$ contains all pairs $(n,k)$ such that $4\leq k\leq n-4$. We also prove that  $4$ is  the least integer $k$ such that every graph $G$ having a large number $n$ of vertices is
$k$-reconstructible up to complementation; this answers a question raised by P. Ille \cite{IL}.\end{abstract}


\section{Introduction}
Ulam Reconstruction Conjecture \cite{Ul} (see \cite{BH, Bo2}) asserts that two graphs $G$ and $G'$ on the same finite set $V$ of $v$ vertices, $v\geq 3$, are isomorphic
provided that the restrictions $G_{\restriction K}$ and  $G'_{\restriction K}$ of $G$ and $G'$ to the
$(v-1)$-element subsets of $V$ are  isomorphic. If this latter condition holds
for the $k$-element subsets of $V$ for some  $k$, $2\leq k \leq v-2$, then, as
it has been noticed several times, $G$ and $G'$ are identical. This conclusion does not requires the finiteness of $v$ nor the isomorphy of  $G_{\restriction K}$ and  $G'_{\restriction K}$, it only requires that 
$G_{\restriction K}$ and  $G'_{\restriction K}$ have the same number of edges for all $k$-element subsets $K$ of $V$,
simply because the adjacency matrix of the Kneser graph $KG(2,k+2)$ is
non-singular (see Section 2).\\

In this paper we look for similar results if the conditions on the restrictions $G_{\restriction K}$
and  $G'_{\restriction K}$ are given up to complementation, that is if $G'_{\restriction K}$ is isomorphic
to  $G_{\restriction K}$ or to its complement ${\overline {G}}_{\restriction K}$, or if $G'_{\restriction K}$ has the same number of edges than $G_{\restriction K}$ or ${\overline {G}}_{\restriction K}$. If the first
condition holds for all $k$-element subsets $K$ of $V$, we say that
$G$ and $G'$ are {\it $k$-hypomorphic up to complementation} and, if the second holds, we say that $G$ and $G'$ have {\it the same number of edges up to complementation}. We say that
$G$ is  {\it $k$-reconstructible  up to complementation}  if every graph $G'$,
$k$-hypomorphic to  $G$ up to complementation, is  isomorphic  to $G$ or its complement. \\

We show first that the equality  of the number of edges, up to complementation, for the $k$-vertices induced subgraphs suffices for the equality up to complementation provided that   $4\leq k\not =7$ and $v$ is  large enough (Theorem \ref{r2(k)}). Our proof is based on Ramsey's theorem for pairs \cite{Ra}.

Next, we give partial description of the set
$\mathcal S$ of pairs $(v,k)$ such that  two graphs $G$ and $G'$ on the same set of $v$ vertices are equal up to complementation whenever they are $k$-hypomorphic up to complementation.

\begin{theorem}\label{equality}
\begin{enumerate}
\item  Let  $v \leq 2$, then $(v,k)\in \mathcal S$ iff $k \in \N$. 
\item Let $v>2$ then $(v,k)\in \mathcal S$ implies $4\leq k\leq v-2$.
\begin{enumerate}
\item If $v\equiv  2 \ (mod \ 4)$,  $(v,k)\in \mathcal S$ iff $4\leq k\leq v-2$; 
\item If $v\equiv  0\ (mod \ 4)$ or $v\equiv  3 \ (mod \ 4)$  then $(v,k)\in \mathcal S$ implies $k\leq v-3$ for infinitely many $v$  and $4\leq k \leq v-3$ implies $(v,k)\in \mathcal S$; 
\item  If $v\equiv  1 \ (mod \ 4)$ then $(v,k)\in \mathcal S$ implies $k\leq v-4$ for infinitely many $v$ and $4\leq k \leq v-4$ implies $(v,k)\in \mathcal S$. \end{enumerate}
 \end{enumerate}
\end{theorem}
Our proof for membership in $\mathcal S$ is  a straithforward application of properties of incidence matrices  due to D.H. Gottlieb \cite{Go}, W. Kantor \cite{KA} and R.M. Wilson \cite{W}. It is given in Section $3$. Constraints on $\mathcal S$ are given in Section $4$.

Our motivation comes from the following problem raised by P. Ille: find the least integer $k$
such that every graph $G$ having a large number $v$ of vertices is
$k$-reconstructible up to complementation. With Theorem \ref{equality} we  show  that
$k=4$ (see Section $2$). 

A quite similar problem was raised by J.G. Hagendorf (1992) and solved by
J.G. Hagendorf and G. Lopez \cite {HL}. Instead of graphs, they consider binary relations and
instead of the complement of a graph, they consider the
{\it dual} $R^*$ of  a binary relation $R$  (where $(x,y)\in R^*$ if and only
if $(y,x)\in R)$; they prove that $12$ is the least integer $k$ such that
 two binary relations $R$ and $R'$, on the same
large set of vertices,  are either isomorphic or dually isomorphic
provided that the  restrictions $R_{\restriction K}$ and $R'_{\restriction K}$ are isomorphic or dually
isomorphic, for every $k$-element subsets $K$ of $V$.
\section{Preliminaries}
Our notations and terminology follow \cite {Bo}.  A {\it graph} is a pair $G:= (V, \mathcal E)$, where $\mathcal E$ is a subset of $[V]^2$, the set of pairs $\{x,y\}$ of distinct elements of $V$. Elements of $V$ are the {\it vertices} of $G$ and elements of $\mathcal E$ its {\it edges}.  If $K$ is a subset of $V$, the {\it  restriction} of $G$ to $K$, also called the {\it induced graph} on $K$ is the graph $G_{\restriction K}:= (K, [K]^2\cap \mathcal E)$. If $K=V\setminus \{x\}$, we denote this  graph by $G_{-x}$. The {\it complement} of  $G$ is the graph  $\overline G:= (V, [V]^2\setminus \mathcal E)$. We denote by $V(G)$ the vertex set of a graph $G$, by $E(G)$ its edge set and by $e(G):=\vert E(G)\vert $ the number of edges. If $\{x,y\}$ is an edge of $G$ we set $G(x,y)=1$; otherwise we set $G(x,y)=0$. The {\it degree} of a vertex $x$ 
of $G$, denoted $d_G(x)$, is  the number of  edges which contain $x$. The graph $G$ is {\it regular} if $d_G(x)=d_G(y)$ for all $x,y\in V$. If $G, G'$ are two graphs, 
we denote by $G\simeq  G'$  the fact that  they are  isomorphic.  A graph is {\it self-complementary} if it is isomorphic to its complement.
\subsection{Incidence matrices and isomorphy up to complementation}
Let $V$ be a finite set, with $v$ elements. Given non-negative integers $t,k$, let $W_{t\;k}$ be the  ${v \choose t}$  by  ${v \choose k}$ matrix of $0$'s and $1$'s, the rows of which are indexed by the $t$-element subsets
 $T$ of $V$, the columns are indexed by the $k$-element subsets $K$
of $V$, and where the entry $W_{t\;k}(T,K)$ is $1$ if
 $T\subseteq K$ and is $0$ otherwise. 

A fundamental result, due to D.H. Gottlieb \cite{Go}, and independently W. Kantor \cite {KA}, is this:  

\begin{theorem}\label{gottlieb-kantor} For $t\leq min{(k,v-k)}$,  $W_{t\; k}$ has full row rank over the field $\Q$ of rational numbers.
 \end{theorem}

 If $k:=v-t$ then, up to a relabelling, $W_{t\; k}$ is the adjacency matrix $A_{t,v}$
 of the  {\it Kneser graph} $KG(t,v)$, graph  whose vertices are the $t$-element
subsets of $V$, two subsets forming an edge if there are disjoint. 

%The adjacency matrix of $KG(k,v)$ is the matrix $A_{k,,k}$ whose rows and columns
%are indexed by the $k$-element subsets of $V$, the coefficient $a_{P,Q}$ being
%$1$ if $P$ and $Q$ are disjoint, $0$ otherwise. We recall a well-known result
%(Gottlieb \cite {Go}, W. Kantor \cite{KA}).
An equivalent form of Theorem \ref{gottlieb-kantor} is:
\begin{theorem}  \label{Ka} $A_{t, v}$ is non-singular for $t\leq \frac{v}{2}$.
\end{theorem}

Applications to graphs and relational structures where given in  \cite{Fr2}  and \cite{Pm}.

 Theorem \ref{gottlieb-kantor} has a  
modular version due to  R.M. Wilson \cite{W}.

\begin{theorem}  For $t\leq min{(k,v-k)}$, the rank of $W_{t\; k}$ modulo a prime $p$ is
$$ \sum  {v  \choose i}-  {v\choose 
i - 1}
 $$
where the sum is extended over those indices $i$ such that $p$
does not divide the binomial coefficient ${k-i \choose  
t-i}$.
\end{theorem}

In the statement of the theorem,  ${v \choose -1}$
should be interpreted as zero.\\

We will apply Wilson's theorem with $t=p=2$ for  $k \equiv 0 \ (mod \ 4)$ and for $k \equiv 1 \ (mod \ 4)$. In the first  case the rank of  $W_{2\; k}$ $(mod \ 2)$ is  
${v  \choose 2} -1$. In the second case, the rank is  ${v \choose  2} -v$.

Let us explain why the use of these results in our context is natural. 

 Let $X_1,\cdots ,X_r$ be an
enumeration of the   $2$-element subsets of $V$; let $K_1,\cdots ,K_{s}$ be an enumeration of the   $k$-element
subsets of $V$ and  $W_{2\; k}$ be the matrix of the $2$-element subsets versus the $k$-element subsets.  If $G$ is a graph with vertex set $V$, let $w_G$ be the row matrix $(g_1,\cdots , g_r)$ where $g_i=1$ if       $X_i$ is an edge of $G$,  $0$ otherwise.   We have
   $w_GW_{2\;Êk}=(e(G_{\restriction K_1}),\cdots , e(G_{\restriction K_s}))$. Thus, if $G$ and $G'$ are two graphs with vertex set $V$ such that  $G_{\restriction K}$ and $G'_{\restriction K}$ have the same number of edges  for every $k$-element subset of $V$, we have  $(w_G-w_{G'})W_{2\;Êk}=0$. Thus, provided that $v\geq 4$, by Theorem \ref{gottlieb-kantor}, $w_G-w_{G'}=0$ that is $G=G'$. 
   
This proves the observation made at the beginning of our introduction.  The same line of proof gives: 

\begin{proposition} \label{down} Let  $t \leq min{(v,  v-k)}$ and $G$ and $G'$ be two graphs on the same set $V$ of $v$ vertices.  If $G$ and $G'$ are $k$-hypomorphic up to complementation then they are $t$-hypomorphic up to complementation.
\end{proposition}
\Proof Let $H$ be a graph on $l$ vertices. Set $Is (H, G):= \{L\subseteq V: G_{\restriction L}\simeq H\}$,  $Isc(H,G):=Is(H,G)\cup Is(\overline H,G)$ and $w_{H,G}$ the $0-1$-row vector indexed by the $t$-element subsets $X_1,\cdots, X_r$ of $V$ whose coefficient of $X_i$ is $1$ if $X_i\in Isc(H,G)$ and $0$ otherwise.  From our  hypothesis, it follows that  $w_{H,G} W_{t \; k}= w_{H,G'} W_{t \; k}$. From Theorem \ref{gottlieb-kantor}, this implies $w_{H,G}= w_{H,G'}$ that is $Isc(H,G)=Isc(H,G')$. Since this equality holds  for all graphs $H$ on $t$-vertices, the conclusion of the lemma follows.
\endproof

Now, let $4\leq k\leq v-4$ and  $G$, $G'$ be two graphs on the same set $V$ of $v$ vertices which are $k$-hypomorphic up to complementation. Then, as shown by Proposition \ref{down}  these two graphs are $4$-hypomorphic up to complementation. By a carefull case analysis (or a very special case of Wilson's theorem, see Theorem \ref{k=0[4]} below), one can prove that two graphs on $6$ vertices  which are $4$-hypomorphic up to complementation are in fact  equal up to complementation. 
Hence, 

\begin{theorem}\label{ille} 
$(k, v)\in \mathcal S$ for all $v, k$ such that $4\leq k\leq v-4$. 
\end{theorem}

P.Ille \cite{IL} asked for the least integer $k$
such that every graph $G$ having a large number $v$ of vertices is
$k$-reconstructible up to complementation. 

From Theorem \ref {ille} above,  $k$ exists and is  at most $4$. From Proposition \ref{nontheo} below, we have $k\geq 4$. Hence $k=4$.

This was our original solution of Ille's problem. The use of Wilson's theorem  leads to the improvement of Theorem \ref{ille} contained in Theorem \ref{equality}. Its use is natural too. Indeed, if we look at conditions which imply $G'=G$ or $G'= \overline G$, it is  simpler to consider the {\it boolean sum}  $G\dot{+} G'$ 
of $G$ and $G'$,  that is the graph $U$ on $V$ whose edges are pairs $e$ of
vertices such that $e\in E(G)$ if and only if $e\notin E(G')$. Indeed,  $G'=G$ or $G'= \overline G$ amounts to the fact that $U$ is either  the empty graph or the complete graph.  But then,  the use of $W_{2\;Êk}$ $(mod \ 2)$ become natural. Particularly,  if conditions which insure 
$w_UW_{2\;Êk}=(0,\cdots ,0)$ $(mod \ 2)$   yield to  $w_U=(0,\cdots ,0)$ or  $w_U=(1,\cdots ,1)$, that is  $U$ is empty or complete, so  $G'=G$ or $G'=\overline G$.

For example,  we show first that if the parity of $e(G_{\restriction K})$ is the same than  $e(G'_{\restriction K})$  for all $k$-element subsets $K$ of $V$, then this may suffice to obtain $G'= G$ or $G'= \overline G$. 

 \begin{theorem}\label{k=0[4]}
Let $G$ and $G'$ be two graphs on the same set $V$ of $v$
vertices (possibly infinite). Let  $k$ be an integer such that $4\leq k\leq v-2$,
$k\equiv 0 \ (mod \ 4)$. Then the following properties are equivalent:\\
(i) $e(G_{\restriction K})$ has the same parity than   $e(G'_{\restriction K})$  for  all $k$-element subsets $K$ of $V$;\\
(ii) $G'= G$ or $G'= \overline G$.
\end{theorem}
\Proof

The implication $(ii)\Rightarrow (i)$ is trivial. We prove  $(i)\Rightarrow (ii)$.\\

\begin{lemma} \label{lem0-4}  If $k \equiv 0  \ (mod \ 4)$ or $k \equiv 1  \ (mod \ 4)$, then under condition (i), $e(U_{\restriction K})\equiv 0  \ (mod \ 2)$ for all $k$-element subsets $K$ of $V$.
\end{lemma}
\Proof
We have trivially :

\begin{claim} \label{congru04} Let $G$ be a graph of $k$ vertices,
then   $e(G)+e(\overline G)$ is even iff  $k \equiv 0  \ (mod \ 4)$ or
 $k \equiv 1  \ (mod \ 4)$.
 \end{claim}

\begin{claim} \label{UG} Let $G$ and $G'$ be two graphs on the same $k$-element vertex set $V$ and let $G\cap G':= (V, E(G)\cap E(G'))$, then :
$$e(G\dot{+} G')=e(G)+e(G')-2e(G\cap G')$$
 \end{claim}

From this, we get :

\begin{claim} \label{ev} Let $G$ and $G'$ be two graphs on a $k$-element  vertex
set $V$, $k\equiv 0  \ (mod \ 4)$ or  $k \equiv 1  \ (mod \ 4)$, and let  $U:=G\dot{+} G'$.   If $e(G')\equiv e(G)  \ (mod \ 2)$ or $e(G')+e(G) \equiv 0  \ (mod \ 2)$ then  $e(U)\equiv 0  \ (mod \ 2)$.
 \end{claim}

The conclusion of Lemma \ref{lem0-4} follows.\endproof

In order to prove implication $(i)\Rightarrow (ii)$ we may suppose $V$ finite. 
 With the notations above, we have $w_UW_{2\;Êk}=(e(U_{\restriction K_1}),\cdots , e(U_{\restriction K_s}))$. Thus,  by Lemma  \ref{lem0-4},  $w_UW_{2\;Êk}=(0,\cdots ,0)$ modulo $2$. Since by Wilson's theorem, the rank of $W_{2\; k}$ modulo $2$ is  ${v\choose 2} -1$,   the kernel of  its transpose $^tW_{2\; k}$ 
 has dimension $1$. Since $(1,\cdots ,1)W_{2\; k}=(0,\cdots ,0)$ $(mod \ 2)$  
  then   $w_UW_{2\; k}=(0,\cdots ,0)$ $(mod \ 2)$ amounts to $w_U=(0,\cdots ,0)$ or  $w_U=(1,\cdots ,1)$, that is  $U$ is empty or complete, so  $G'=G$ or $G'=\overline G$.
\endproof\\
 \begin{remark}  For every integer $k \not\equiv 0  \ (mod \ 4)$ there are two graphs
 $G$ and $G'$ on the same vertex set $V$, $\vert V\vert \geq k+2$,  such that
 $e(G_{\restriction K})$ has the same parity than   $e(G'_{\restriction K})$  or ${{k(k-1)}\over 2}- e(G'_{\restriction K})$ for all
$k$-element subsets $K$ of $V$, but  $G$ and $G'$ are not
isomorphic  up to complementation.
\end{remark}

For an example, consider $G$ and $G'$ on the same vertex set $V:=\{1,\cdots ,v\}$ such that the edges of $G'$ are $(1,i)$  for all $i\in  \{2,3,\cdots ,v\}$ and $G$ is the empty graph if $k\equiv 3  \ (mod \ 4)$ or $k\equiv 1  \ (mod \ 4)$; $G$ is a complete graph if $k\equiv 2  \ (mod \ 4)$.\\

 We give  an analog of Theorem \ref{k=0[4]} in the case $k\equiv 1 \ (mod \ 4)$. For that, an additional condition  is needed. 

Let $G$ be a graph.  A $3$-element subset  $T$ of $V$ such that all
pairs belong to $E(G)$ is a 
 {\it triangle} of $G$. A $3$-element subset of $V$ which is a triangle of $G$ or of $\overline G$ is a  $3$-{\it homogeneous} subset of $G$.

 \begin{theorem}\label{k=1[4]}
Let $G$ and $G'$ be two graphs on the same set $V$ of $v$
vertices (possibly infinite). Let  $k$ be an integer such that $5\leq k\leq v-2$,
$k\equiv 1 \ (mod \ 4)$. Then the following properties are equivalent:\\
(i) $e(G_{\restriction K})$ has the same parity than   $e(G'_{\restriction K})$  for  all $k$-element subsets $K$ of $V$ and the same $3$-homogeneous subsets;  \\
(ii) $G'= G$ or $G'= \overline G$.
\end{theorem}

\Proof The implication $(ii)\Rightarrow (i)$ is trivial. We prove  $(i)\Rightarrow (ii)$.\\
We may suppose $V$ finite. Let $U:= G\dot{+} G'$. From the fact that $e(G_{\restriction K})$ has the same parity than   $e(G'_{\restriction K})$  for  all $k$-element subsets $K$, the boolean sum $U$ belongs to the kernel of  $^tW_{2\; k}$  (over the $2$-element field). 

\begin{claim} \label{bipartite}Let $k$ be an integer such that $2\leq k\leq v-2$,  $k\equiv 1 \ (mod \ 4)$,  then the kernel of $^t W_{2\; k}$ consists of complete bipartite graphs and their complements (including the empty graph and the complete graph).
\end{claim}
\Proof  Let us recall that a {\it star-graph} of $v$ vertices consists of a vertex linked to all other vertices, those $v-1$ vertices forming an independent set. The vector space (over the $2$-element field) generated by the star-graphs on $V$ consists of all complete bipartite graphs distinct from the complete graph (but including  the empty graph). Moreover,  its dimension is $v-1$ (a basis being made of star-graphs).  Let $\K$ be the kernel of $^t W_{2\; k}$.  Since $k$ is odd, each star-graph belongs to $\K$. Since $k\equiv 1 \ (mod \ 4)$, the complete graph also belongs to $\K$. According to Wilson's theorem,  the rank of  $W_{2k}$ $(mod \ 2)$ is  
$ {v\choose 2} -v$. Hence  the kernel of $^tW_{2\; k}$ has dimension $v$.  Consequently, $\K$ consists of complete bipartite graphs and their complements, as claimed. 
 \endproof
 
A {\it claw} is a star-graph on  four  vertices, that is a graph made of a vertex joined to three other vertices, with no edges between these  three vertices.  A graph is {\it claw-free} if no induced subgraph is a claw. 

\begin{claim} \label{clawfree}Let $G$ and $G'$ be two graphs on the same set and having the same $3$-homogeneous subsets, then the boolean sum $U: =G\dot {+} G'$ and its complement are  claw-free. 
\end{claim}
\Proof Let $x\in V$. Suppose $d_U(x)\geq 3$. Then, the  neighborhood of  $x$ contains  at least two distinct vertices $y,y'$ such that $U(y,y')=1$. Indeed, it contains clearly two vertices $y,y'$ such that  $G(x,y)=G(x,y')$. If  $U(y,y')=0$,  that is  $G(y,y')= G'(y,y')$,  then since $G$ and $G'$ have the same $3$-element homogeneous sets and $G(x,y)\not =G'(x, y)$, $\{x,y,y'\}$ cannot be homogeneous, hence $G(y,y')\not = G(x,y)$ and $G'(y,y')\not = G'(x,y)$. This implies $G(y,y')\not =G'(y,y')$,  a contradiction. 
 From this observation, $U$ is claw-free. Since $G$ and $\overline G'$ have the same $3$-homogeneous subsets and $\overline U= G\dot{+}\overline G'$, we also get  that $\overline U$ is claw-free.  
\endproof

For a characterization of these boolean sums, see \cite {PS}.

From Claim \ref{bipartite}, $U$ or its complement is a complete bipartite graph and,  from Claim \ref{clawfree},  $U$ and $\overline U$ are claw-free. Since $v\geq 5$ (in fact $v\geq 7$), it follows that $U$ is either the empty graph or the complete graph. Hence $G'=G$ or $G'= \overline {G}$ as claimed. 
\endproof

\subsection{Conditions on the number of edges and Ramsey's theorem}

\begin{theorem}\label{r2(k)}
Let  $k$ be an integer,  $7\neq k\geq 4$. There is an integer $m$ such that 
if $G$ and $G'$ are two graphs on the same set $V$ of $v$ vertices, $v\geq m$,   such that
$G_{\restriction K}$ and $G'_{\restriction K}$  have the same number of edges, up to
complementation, for  all $k$-element subsets $K$ of $V$, then $G'= G$ or $G'= \overline G$.
\end{theorem}

Conditions $7\neq k\geq 4$ in Theorem \ref{r2(k)} are  necessary.

- For $k=7$, consider two graphs $G$ and $G'$ on
$V:=\{1,2,\cdots ,v\}$   such that $\{i,j\}$ is an edge of $G$ and $G'$ for all $i\neq j$ in $\{1,2,\cdots ,v-2\}$,
$G$ has no another edge and $G'$ has $\{v-1,v\}$ as an additional edge. For $k\geq 4$ apply Proposition \ref{nontheo} below.

Let  $c(k)$ be  the least integer $m$ for which the conclusion of  Theorem \ref{r2(k)}  holds. 

\begin{problem}  Is $c(k) \leq k+4$?
\end{problem} 
 Our proof uses Ramsey's theorem rather than incidence matrices.  It is inspired from a relationship between Ramsey's theorem  and Theorem \ref{gottlieb-kantor} pointed out in \cite{Pm}.  The drawback is that the bound on $c(k)$ is quite crude.

Let $r_2^2 (k)$ be the bicolor Ramsey number  for pairs: the least integer $n$ such that every graph on $n$ vertices contains a $k$-homogeneous subset, that is a clique or an independent on $k$ vertices.  We  deduce Theorem \ref{r2(k)} and 
$c(k)\leq  r_2^2 (k)$ from the  following result.
\begin{proposition}\label{prop{r2(k)}}
Let  $k$ be an integer,  $7\neq k\geq 4$ and let $G$ and $G'$ be  two graphs on the same set $V$ of $v$ vertices, $v\geq k$ such that:
\begin{enumerate}
 \item$G_{\restriction K}$ and  $G'_{\restriction K}$ have the same number of edges, up to complementation, for  all $k$-element subsets $K$ of $V$;
\item $V$ contains a $k$-element subset $K$ such that $G_{\restriction K}$ or  ${\overline G}_{\restriction K}$ has at least $l$ edges where $l:= \min \left( {{k^2+7k-12} \over 4}, {{k(k-1)}\over 2}\right)$.
\end{enumerate}
Then $G'= G$ or $G'= \overline G$.
\end{proposition}

The inequality $ {{k^2+7k-12} \over 4}\leq {{k(k-1)}\over 2}$ holds iff  $k\geq 8$. For $k> 8$ the condition $l= {{k^2+7k-12} \over 4}$ is weaker than the existence of a clique of size $k$.
 
\Proof
We may suppose that $V$ contains a  $k$-element subset   of $V$, say  $K$, such that
 $e(G_{\restriction K}) \geq l$; also we may suppose, 
 from condition 1, that $e(G_{\restriction K})=e(G'_{\restriction K})$  
 otherwise replace  $G'$ by its complement. 
We shall prove that for all $V'$ such that $K\subseteq V'\subseteq V$ and $\mid V'\mid = k+2$ we have
$e(G_{\restriction K'})=e(G'_{\restriction K'})$ for all $k$-element subset $K'$ of $V'$. Since the adjacency matrix of the Kneser graph $KG(2,k+2)$ is
non-singular,  $G_{\restriction   V'}=G'_{\restriction   V'}$. It follows that $G=G'$. 

\begin{claim} \label{c1} For $x\notin K$ and $y\in K$,  $e(G_{\restriction   (K\cup\{x\})\setminus \{y\}}) = e(G'_{\restriction   (K\cup\{x\})\setminus \{y\}})$. \end{claim}
\Proof   Let    $x\notin K$ and $y\in K$. Set $K':=(K\cup\{x\})\setminus \{y\}$.  The graphs $G_{\restriction   K'}$ and
$G'_{\restriction   K'}$ have at least $l':=l - (k-1)$ edges. Since  $G_{\restriction K'}$ and $G'_{\restriction K'}$ have the same number of edges up to complementation, we have $e(G_{\restriction  K'})= e(G'_{\restriction  K'})$ whenever $l'\geq  {{k(k-1)}\over 4}$, that is $l\geq l'':= {{(k-1)(k+4)}\over 4}$.

If $k\geq 8$ we have  $ l= {{k^2+7k-12}\over 4}$ yielding $l>l''$ as required. If $k\in \{4,5,6\}$ we have $l= {{k(k-1)}\over 2}$ yielding again $l\geq l''$. \endproof

\begin{claim} \label{c1bis} For distinct $x,x'\notin K$ and $y,y'\in K$,  $e(G_{\restriction   (K\cup\{x,x'\})\setminus \{y,y'\}}) = e(G'_{\restriction   (K\cup\{x,x'\})\setminus \{y,y'\}})$. \end{claim}
\Proof 
Let $x,x'\notin K$ and $y, y'\in K$ be distinct. Set $K':=(K\cup \{x,x'\}) \setminus
\{y,y'\}$. We have $e(G_{\restriction K'})\geq e(G_{\restriction K})-(2k-3)$ and $e(G'_{\restriction K'})\geq
e(G_{\restriction K})-(2k-3)$. Thus  $e(G_{\restriction K'})$ and   $e(G'_{\restriction K'})$ have at least
$l' := l-(2k-3)$ edges.  Since  $G_{\restriction K'}$ and $G'_{\restriction K'}$ have the same number of edges up to complementation, we have $e(G_{\restriction  K'})= e(G'_{\restriction  K'})$ whenever $l'\geq  {{k(k-1)}\over 4}$, that is $l\geq  {{k^2+7k-12} \over 4}$. This inequality holds if $k\geq 8$.

Suppose $k\in \{4,5,6\}$. Thus $l= {{k(k-1)}\over 2}$. Hence $K$ is a clique for $G$ and $G'$.  

\noindent {\bf Subclaim} Let $u\notin K$ then $G$
 and $G'$ coincide on $K\cup \{u\}$. 
 
 \Proof Since $K$ is a clique, this amounts to $G(u,v) = G'(u,v)$  for all $v\in K$, a fact which follows from Claim \ref{c1}. Indeed, we have $d_{G_{\restriction K\cup\{u\}}}(u)= {1\over {k-1}}\sum_{w\in K} d_{G_{\restriction (K\cup\{u\})\setminus \{w\}}}(u)$. From Claim \ref{c1} we have 
$d_{G_{\restriction  (K\cup\{u\})\setminus \{w\}}}(u)=d_{G'_{\restriction  (K\cup\{u\})\setminus \{w\}}}(u)$. Thus $d_{G_{\restriction K\cup\{u\}}}(u)=d_{G'_{\restriction K\cup\{u\}}}(u)$. Since $d_{G_{\restriction  (K\cup\{u\})\setminus \{v\}}}(u)=d_{G'_{\restriction  (K\cup\{u\})\setminus \{v\}}}(u)$ the equality  $G(u,v) = G'(u,v)$ follows.\endproof

 From this subclaim it follows that $G$ and $G'$ coincide on $K'$ with the possible exception of the pair $\{x,x'\}$. Set  $a:=e(G_{\restriction   K'})$, 
$a':=e(G'_{\restriction   K'})$. Suppose  $a\neq a'$.  Then $\vert a-a'\vert =1$, hence  the sum $a+a'$ is odd. Since 
$G_{\restriction   K'}$ and $G'_{\restriction   K'}$ have the same number of edges up to complementation, this sum is also ${{k(k-1)}\over 2}$.  If $k=4$ or $k=5$ this number is even, a contradiction. 
 Suppose $k=6$. We may suppose $a=a'+1$ hence from $a+a'={{k(k-1)}\over 2}$ we get $a=8$ .
Put $\{x_1,x_2,x_3,x_4,y,y'\}:=K$. Since $K$ is a clique we have $G(x,x')=1$, $G'(x,x')=0$ and $G$, $G'$ contain just an edge from $\{x, x'\}$ to $\{x_1,x_2,x_3,x_4\}$. We  may suppose $G(x_1,x)=G'(x_1,x)=1$, $G(x_1,x')=G'(x_1,x')=0$ and $G(t,u)=G'(t,u)=0$ for all $t\in \{x_2,x_3,x_4\}$ and $u\in \{x,x'\}$. \\
Let $K'':=(K\cup \{x,x'\})\setminus \{x_1,x_2\}$. From the subclaim above, $G$ and $G'$ coincide on $K''$ at the exception of the pair $\{x,x'\}$ hence $G$, $G'$ contain   just an edge from $\{x, x'\}$ to $\{x_3,x_4, y,y'\}$. We can assume $G(y,u)=G'(y,u)=1$ for exactly one $u\in \{x,x'\}$, and  $G(t,u)=G'(t,u)=0$ for all $t\in \{x_3,x_4,y'\}$ and $u\in \{x,x'\}$.

Set $B:=\{x_2,x_3,x_4,x,x',y'\}$, then $e(G_{\restriction B})=7$ and
$e(G'_{\restriction B})=6$. So $e(G_{\restriction B})\neq e(G'_{\restriction B})$ and 
$e(G_{\restriction B})+e(G'_{\restriction B})\neq  {{k(k-1)}\over 2}$, that gives a contradiction.\endproof       

\noindent Clearly  Proposition \ref{prop{r2(k)}} follows from Claims  \ref{c1}   and  \ref{c1bis}.
 \endproof

\section{Some members of $\mathcal S$}

Sufficient conditions for membership stated in Theorem \ref{equality} are contained in Theorem \ref{principal} below.

Let $v$ be a non negative integer and $\vartheta(v):= 4l$ if $v\in\{ 4l+2, 4l+3\}$, $\vartheta (v):=4l-3$ if $n\in\{4l, 4l+1\}$.
\begin{theorem}\label{principal}
Let $v$, $k$ be  two integers , $v\geq 6$, $4\leq k \leq \vartheta (v)$. Then, for every pair of graphs  $G$ and $G'$ on the the same set $V$ of $v$ vertices, the following properties are equivalent:
\begin{enumerate}[{(i)}]
\item $G$ and  $G'$ are $k$-hypomorphic up to complementation;
\item $G_{\restriction K}$ and  $G'_{\restriction K}$ have the same number of edges, up to complementation,
and the same number of $3$-homogeneous subsets, for  all $k$-element
subsets $K$ of $V$;
\item  $G_{\restriction K}$ and  $G'_{\restriction K}$ have the same number of edges, up to complementation, for
all $k$-element and $k'$-element subsets $K$ of $V$ where $k'$ is an integer verifying $3\leq k'<k$;
\item $G'= G$ or $G'= \overline G$.
\end{enumerate}
\end{theorem}

\subsection{Ingredients}
Let $G:=(V,E)$ be a graph.   Let
$A^{(2)}(G)$ be the set  of
unordered pairs $\{u,u'\}$ made of some $u \in E(G)$ and some $u'\in
E(\overline G)$.  Let   $A^{(0)}(G):=\{ \{u,u'\} \in A^{(2)}(G)\ : \ u\cap
u'=\emptyset \}$,  $A^{(1)}(G):=A^{(2)}(G)\setminus A^{(0)}(G)$ and let $a^{(i)}(G)$ be
the cardinality of  $A^{(i)}(G)$ for $i\in \{0,1,2\}$; thus  $a^{(2)}(G)= a^{(0)}(G)+a^{(1)}(G)$.
Let $T(G)$ be the set
of {\it triangles} of $G$ and let $t(G):=\mid T(G)\mid$. Let
$H^{(3)}(G):=T(G)\cup T(\overline G)$ be the set of $3$-homogeneous subsets of $G$ and $h^{(3)}(G):=\mid H^{(3)}(G)\mid$.
%The relevance of these notions is due  to the following  facts.
% \begin{fact} \label{trivial}
%Let $G$ and $G'$ be two graphs on a $3$-element vertex set.Then the following properties are equivalent:
%\begin{enumerate}[{(i)}]
%\item $G$ and $G'$ are $3$-hypomorphic up to complementation;
%\item $G$ and $G'$ have the same number of edges up to complementation;
%\item  $H^{(3)}(G)=H^{(3)}(G')$.  
%\end{enumerate}
%\end{fact}
%\begin{fact} \label{lem000}

Some elementary properties of the above numbers are stated in the lemma below;  the proof is  immediate. 

\begin{lemma} \label{claim} Let $G$ be a graph with $v$ vertices, then :\\
1) $A^{(i)}(G)=A^{(i)}(\overline G)$, hence
$a^{(i)}(G)=a^{(i)}(\overline G)$, for all $i\in \{0,1,2\}$.\\
2) $a^{(2)}(G)=e(G)e({\overline G})$.\\
3) $a^{(1)}(G)=\sum_{x\in V(G)}d_G(x)d_{\overline G}(x)$.\\
4) $h^{(3)}(G)= {{v(v-1)(v-2)}\over {6}}-{1\over 2}a^{(1)}(G)$.
\end{lemma}

\begin{lemma} \label{lem00}
Let $G$ and $G'$ be two graphs on the same finite vertex set $V$, then  :
 $$e(G')= e(G)\ or \ e(G')=e({\overline G})  \  \mbox{iff}  \
 e(G)e({\overline G})=e(G')e({\overline {G'}})$$
\end{lemma}

\Proof Suppose :
\begin{equation}\label{eq1}
e(G)e({\overline G})=e(G')e({\overline {G'}})
\end{equation}
Since $e(G)+e({\overline G})= {{v(v-1)}\over 2}$ and $e(G')+e({\overline
{G'}})= {{v(v-1)}\over 2}$, where $v:=\mid V\mid$, we have :
\begin{equation}\label{eq2}
e(G)+e({\overline G})=e(G')+e({\overline {G'}}) 
\end{equation}
Then (\ref{eq1})  and  (\ref{eq2}) give $e(G')= e(G)\ or \ e(G')=e({\overline G})$. The converse is obvious. 
\endproof

\begin{lemma} \label{lem0000}
Let $G$ be a graph, $V:=V(G), v:=\mid V\mid$.\\
a) Let $i\in \{ 0,1\}$, $k$ such that $4-i\leq k\leq v$, then :
$$a^{(i)}(G)= {1\over {v-4+i \choose k-4+i}}\sum_{\s K\subseteq V \atop  \mid K\mid =k}
a^{(i)}(G_{\restriction K})$$
b) Let $k$ such that $3\leq k \leq v-1$, then :
$$a^{(0)}(G)= {{v-3}\over {v-k}}e(G)e({\overline {G}})-{{1}\over
{{v-4\choose k-3}}}
\sum_{ \s K\subseteq V \atop \mid K\mid =k} e(G_{\restriction K})e({\overline G}_{\restriction K})$$
$$a^{(1)}(G)= {{1}\over {{v-4\choose k-3}}}\sum_{\s K\subseteq V \atop \mid K\mid =k}
e(G_{\restriction K})e({\overline G}_{\restriction K})-{{k-3}\over {v-k}}e(G)e({\overline {G}})$$
 \end{lemma}

\Proof
a)  Let $\{u,u'\} \in A^{(i)}(G)$ for $i\in \{0,1\}$. The number of $k$-element subsets $K$  of $V$   containing $u$
and $u'$ is ${v-4+i \choose k-4+i} $. Then we conclude.

b) If $k=3$ then a) and the fact that $a^{(0)}(G)+a^{(1)}(G)=e(G)e({\overline {G}})$ give the formulas.\\
If $4\leq k \leq v-1$, then by a)  we have :
$${{v-4\choose k-4}}a^{(0)}(G)=\sum_{ \s
K\subseteq V \atop  \mid K\mid =k}a^{(0)}(G_{\restriction K})$$
$${{v-3\choose k-3}}a^{(1)}(G)=\sum_{  \s
K\subseteq V \atop \mid K\mid =k}a^{(1)}(G_{\restriction K})$$
Summing up and applying  2) of Lemma \ref{claim} to the $G_{\restriction K}$'s we have :
\begin{equation}\label{eq3}
{{v-4\choose k-4}}a^{(0)}(G) + {{v-3\choose k-3}}a^{(1)}(G)=
\sum_{\s K\subseteq V \atop  \mid K\mid =k} e(G_{\restriction K})e({\overline G}_{\restriction K})
\end{equation}
On an other hand :
\begin{equation}\label{eq4}
a^{(0)}(G)+a^{(1)}(G)=e(G)e({\overline G})
\end{equation}
Equations (\ref{eq3}) and  (\ref{eq4})  form a Cramer system with $a^{(0)}(G)$ and $a^{(1)}(G)$
as
unknowns.  Indeed the determinant
$\Delta := \left|\begin{array}{cc}
{v-4\choose k-4} & {v-3\choose k-3}\\
1 & 1
\end{array}
 \right|
= {v-4\choose k-4} - {v-3\choose k-3}= - {v-4\choose k-3}$ is non zero.
A straightforward computation gives the result. \endproof

\begin{corollary} \label{corkk+1}
Let $G$ and $G'$ be two graphs on the same set $V$ of $v$ vertices and 
$k$ be an integer such that $4\leq k\leq v$. 

 The implications $(ii)\Rightarrow (i)$ and  $(i)\Rightarrow (iii)$ between the following statements hold.
\begin{enumerate}[{(i)}]

\item $e(G'_{\restriction K})= e(G_{\restriction K})$ or $e({\overline G_{\restriction K}})$ and
$h^{(3)}(G_{\restriction K})=h^{(3)}(G'_{\restriction K})$ for all $k$-element subsets $K$ of $V$.

\item  $e(G'_{\restriction K})= e(G_{\restriction K})$ or $e({\overline G_{\restriction K}})$
for all $k$-element and $k'$-element subsets $K$ of $V$ where $k'$ is some  integer  verifying $3\leq k'<k$.

\item $G_{\restriction L}$ and $G'_{\restriction L}$ have the same number of edges up to complementation and
$h^{(3)}(G_{\restriction L})=h^{(3)}(G'_{\restriction L})$ for  all $l$-element  subsets  $L$ of $V$  and all integer $l$ 
such that $k\leq l \leq v$.
\end{enumerate}

 \end{corollary}

\Proof $(i)\Rightarrow (iii)$.  Let $L$  be an $l$-element  subset of $V$ with $l\geq k$, and
 $K$ be a $k$-element subset of $L$.
From Lemma \ref{lem00} and 2) of Lemma \ref{claim}, we have
$a^{(0)}(G_{\restriction K})+a^{(1)}(G_{\restriction K})=a^{(0)}(G'_{\restriction K})+a^{(1)}(G'_{\restriction K})$, and from 4)
 of Lemma \ref{claim},
$a^{(1)}(G_{\restriction K})=a^{(1)}(G'_{\restriction K})$. Hence $a^{(i)}(G_{\restriction K})=a^{(i)}(G'_{\restriction K})$ for all
$k$-element subsets $K$ of $L$ and $i\in \{0,1\}$.\\
From a) of Lemma \ref{lem0000} applied to $G_{\restriction L}$  follows $a^{(i)}(G_{\restriction L})=a^{(i)}(G'_{\restriction L})\ for \ i\in \{
0,1\}$, hence using   2) of Lemma \ref{claim} we get
$e(G_{\restriction L})e({\overline G_{\restriction L}})=e(G'_{\restriction L})e({\overline {G'}_{\restriction L}})$.
The conclusion follows from Lemma  \ref{lem00} and 4) of Lemma  \ref{claim}.\\

$(ii)\Rightarrow (i)$.  It suffices to prove that $h^{(3)}(G_{\restriction K})=h^{(3)}(G'_{\restriction K})$ for all $k$-element subsets $K$ of
$V$. From Lemma \ref{lem00} we have $e(G_{\restriction K})e({\overline
G_{\restriction K}})=e(G'_{\restriction K})e({\overline G'}_{\restriction K})$ and $e(G_{\restriction K'})e(\overline
G_{\restriction K'})=e(G'_{\restriction K'})e({\overline G'}_{\restriction K'})$ for all $k'$-element set
$K'\subseteq K$.  From b) of Lemma \ref{lem0000} we get
$a^{(i)}(G_{\restriction K})=a^{(i)}(G'_{\restriction K})$ for $\ i\in \{ 0,1\}$. Then by 4) of
Lemma  \ref{claim}, $h^{(3)}(G_{\restriction K})=h^{(3)}(G'_{\restriction K})$.  \endproof

\begin{proposition} \label{e(G)=e(G')}
Let $G$ and $G'$ be two graphs on $v$ vertices and $k$ be an integer such that $4\leq k\leq v$. If $G$ and $G'$ are 
$k$-hypomorphic up to complementation
then $e(G'_{\restriction L})= e(G_{\restriction L})\ or \ e(G'_{\restriction L})=e({\overline G}_{\restriction L})$ for all $l$-element  subsets  $L$ of $V$  and all integer $l$ such that 
$ k\leq l\leq v$.
\end{proposition}

\Proof
If $G$ and $G'$ are $k$-hypomorphic up to complementation  then
$G_{\restriction K}$ and  $G'_{\restriction K}$ have the same number of edges up to complementation,
and the same number of $3$-homogeneous subsets, for  all $k$-element
subsets $K$ of $V$. We conclude using $(i)\Rightarrow (iii)$  of Corollary \ref{corkk+1}
\endproof \\
By inspection of the eleven graphs on four vertices, one may observe that:
\begin{fact} \label{lem000}The pair $(e(G)e({\overline G}),h^{(3)}(G))$ characterize $G$ up to isomorphy
and complementation if $\mid V(G)\mid \leq 4$.
 \end{fact}
Note that in  Fact  \ref{lem000}, we can replace $(e(G)e({\overline G}),h^{(3)}(G))$ by  $(a^{(0)}(G),a^{(1)}(G))$ (this follows from Lemmas \ref{lem00} and  \ref{claim}).\\

\begin{proposition}  \label{Ka+} Let $G$ and $G'$ be two graphs on the
same set $V$ of $v$ vertices and $k$ be an integer.
If $3\leq k\leq v-3$ (resp. $4\leq k\leq v-4$) and
$h^{(3)}(G_{\restriction K})=h^{(3)}(G'_{\restriction K})$ (resp.  $a^{(0)}(G_{\restriction K})=a^{(0)}(
{G'}_{\restriction K})$) for all $k$-element subsets $K$ of $V$ then
$h^{(3)}(G_{\restriction K})=h^{(3)}(G'_{\restriction K})$
(resp. $a^{(0)}(G_{\restriction K})=a^{(0)}({G'}_{\restriction K})$)
for all $(v-k)$-element subsets $K$ of $V$.
\end{proposition}

\Proof  By $4)$ of Lemma \ref{claim},  $h^{(3)}(G_{\restriction K})=h^{(3)}(G'_{\restriction K})$ iff
 $a^{(1)}(G_{\restriction K}) =a^{(1)}(G'_{\restriction K})$.\\

Case 1.  $ k\leq {v\over 2}$,      then $v-k\geq k$. Let   $K'$ be a $(v-k)$-element subset of  $V$, then
from a) of Lemma   \ref{lem0000}      we have for  $\ i\in \{ 0,1\}$ :
 $$a^{(i)}(G_{\restriction K'})= {1\over {{v-k-4+i \choose k-4+i}}}\sum_{\s K\subseteq K' \atop \mid K\mid =k}
a^{(i)}(G_{\restriction K})$$
Then we get the conclusion.\\

Case 2.   $ k> \frac {v}{2}$,     then $v-k<   \frac {v}{2}$.
 Let   $K'$ be a $k$-element subset of  $V$.
From a) of Lemma   \ref{lem0000}      we have for  $\ i\in \{ 0,1\}$ :
\begin{equation}\label{eqkneser}\sum_{\s K\subseteq K' \atop \mid K\mid =v-k} a^{(i)}(G_{\restriction K})={k-4+i \choose v-k-4+i} a^{(i)}(G_{\restriction K'})
\end{equation}
 Let $X_1,X_2,\cdots ,X_l$ be an enumeration of the   $(v-k)$-element subsets of  $V$. 
Let $w_G^{(i)}:= (a^{(i)}(G_{\restriction X_1}), a^{(i)}(G_{\restriction X_2}),   \cdots , a^{(i)}(G_{\restriction X_l}))$, and 
$w_{G'}^{(i)}:= (a^{(i)}(G'_{\restriction X_1}), a^{(i)}(G'_{\restriction X_2}),   \cdots , a^{(i)}(G'_{\restriction X_l}))$. 
 From (\ref{eqkneser}), we
 get, for $i\in \{0,1\}$ :
 $A_{v,v-k} {^t w_G^{(i)}}= A_{v,v-k}{^t w_{G'}^{(i)}}$. 
We conclude using Theorem \ref{Ka}.  \endproof

\subsection {Proof of  Theorem \ref{principal}.}

$(i) \Rightarrow (ii)$, $(iv) \Rightarrow (i)$, $(iv) \Rightarrow (iii)$ are obvious and
$(iii) \Rightarrow (ii)$ is implication $(ii)\Rightarrow (i)$ of  Corollary \ref{corkk+1}. Thus  it is sufficient to prove
$(ii) \Rightarrow (iv)$. \\
 Let $l$, $k\leq l \leq v$. According to implication $(i)\Rightarrow (ii)$ of Corollary
\ref{corkk+1}, $e(G'_{\restriction L})=e(G_{\restriction L})$ or $e(G'_{\restriction L})=e(\overline G_{\restriction L})$
for all $l$-element subsets $L$ of $V$. 
If  we may choose  $l\equiv 0  \ (mod \ 4)$ with $l\leq v-2$, then from Claim \ref{congru04},
 $e(G_{\restriction L})$ and $e(G'_{\restriction L})$ have the same parity. Theorem \ref{k=0[4]} gives $G'=G$ or $G'=\overline G$.  Thus,  the implication $(ii)\Rightarrow (iv)$ is proved if   $v\equiv 2 \ (mod \ 4)$ and if $v\equiv 3 \ (mod \ 4)$. There are two remaining cases. 
  
{\bf Case 1. }   $v\equiv 1 \ (mod \ 4)$ and $k=v-4$. We prove that $e(G'_{\restriction L})$ and $e(G_{\restriction L})$ have the same parity for all $4$-element subsets $L$ of $V$. Theorem \ref{k=0[4]} again gives $G'=G$ or $G'=\overline G$. 
The proof goes as follows.  Let $L$  be a $4$-element  subset of $V$, and
 $K$ be a $k$-element subset of $V$.  By  Lemma \ref{claim},  $a^{(2)}(G_{\restriction K})=a^{(2)}(G'_{\restriction K})$
  and   $a^{(1)}(G_{\restriction K})=a^{(1)}(G'_{\restriction K})$. Thus $a^{(0)}(G_{\restriction K})=a^{(0)}(G'_{\restriction K})$. Using Proposition \ref{Ka+},
  we get
  $a^{(0)}(G_{\restriction L})=a^{(0)}(G'_{\restriction L})$ and $h^{(3)}(G_{\restriction L})=h^{(3)}(G'_{\restriction L})$.
  Now 4) of Lemma  \ref{claim} gives $a^{(1)}(G_{\restriction L})=a^{(1)}(G'_{\restriction L})$.
  So $a^{(2)}(G_{\restriction L})=a^{(2)}(G'_{\restriction L})$, then using 2) of Lemma  \ref{claim}
  and Lemma  \ref{lem00} we get $e(G'_{\restriction L})=e(G_{\restriction L})$ or $e({\overline {G}}_{\restriction L})$, thus
   $e(G'_{\restriction L})$ and $e(G_{\restriction L})$ have the same parity.\endproof
   
{\bf Case 2.} $v\equiv 0 \ (mod \ 4)$ and $k=v-3$. From Proposition  \ref{Ka+},  $G$  and $G'$ have the same $3$-homogeneous subsets. From Theorem \ref{k=1[4]},
 $G'=G$ or $G'= \overline {G}$ as claimed. 
\endproof

\section{Constraints on $\mathcal S$ }
Two arbitrary graphs on the same set of vertices are $k$-hypomorphic up to complementation  for $k\leq 2$. Hence, if $v\leq 2$, $(v,k)\in \mathcal S$ iff $k\in \N$. This is item $1$ of Theorem \ref{equality}. 

Next,  suppose  $v> 2$,  and $(v,k)\in \mathcal S$.

According to the proposition below, we have $k\geq 4$. 
\begin{proposition} \label{nontheo} For every integer $v\geq 4$, there are two
graphs $G$ and $G'$, on the same set of $v$ vertices, which are $3$-hypomorphic up to complementation but not isomorphic up to complementation.
\end{proposition}

\Proof
Let  $G$ and $G'$ be two graphs
having $\{1,2,\cdots ,v\}$ as set of vertices.

 \noindent - Even case : $v=2p$. Pairs  $\{i,j\}$  are edges of
$G$ and $G'$ for all $i\neq j$ in $\{1,2,\cdots ,p\}$ and for all $i\neq j$ in
$\{ p+1,\cdots ,2p\}$. The graph $G$ has no other edge and $G'$ has $\{1,
p+1\}$ as an additional edge. Clearly $G'$ and $G$ are $3$-hypomorphic up to complementation and not isomorphic. Since 
$\overline G$ has  $p^2$ edges but $G'$ has $p(p-1)+1$ edges, $G'$ and $\overline G$ are not isomorphic.\\
- Odd case : $v=2p+1$. Pairs $\{i,j\}$ are edges of  $G$ and
$G'$ for all $i\neq j$ in $\{1,2,\cdots ,p\}$ and for all $i\neq j$ in $\{p+1,\cdots ,2p+1\}$.  The graph $G$ has no other edge and $G'$ has $\{1,
p+1\}$ as an additional edge. Clearly $G'$ and $G$ are $3$-hypomorphic up to complementation and not isomorphic.  Since 
$\overline G$ has  $p(p+1)$ edges but $G'$ has $p^2+1$ edges,  $G'$ and $\overline G$ are not isomorphic.\\
In both cases $G$ and $G'$  are $3$-hypomorphic up to complementation but not
isomorphic up to complementation. \endproof

According to the following lemma,  $v\geq 6$. 
\begin{lemma} \label{five}For every $v$, $3\leq v\leq 5$, there are two graphs  $G$ and $G'$,  on the same set of $v$ vertices,   which are $k$-hypomorphic for  all $k\leq v$ but $G'\neq G$ and $G'\neq {\overline G}$.
\end{lemma}\label{five}
\Proof
Let $V:= \{0, 1, 2, 3,4\}$, $\mathcal E:= \{\{0,1\}, \{1,2\}, \{2,3\}, \{3,4\}, \{4,0\}\}$, $\mathcal E':= (\mathcal E \setminus \{\{0,4\},$ $\{1,2\}\})Ê\cup \{\{1,4\}, \{0,2\}\}$.  Let $G:= (V, \mathcal E)$ and  $G':= (V, \mathcal E')$. These  graphs are  two $5$-element cycles, $G'$ being obtained from $G$ by exchanging $0$ and $1$. Trivially, they satisfy the conclusion of the lemma. The two pairs  $G_{-3}$ ,  $G'_{-3}$ and $G_{-3, -4}$ and $G'_{-3,-4}$ also satisfy the conclusion of the lemma.
\endproof

Next, a  straightforward extension  of the construction in Lemma \ref{five} above   yields $k\leq v-2$. 
Indeed, let us say that  two graphs $G$ and $G'$ on the same set $V$ of vertices are {\it $k$-hypomorphic} if for any subset $X$ of
$V$ of cardinality $k$,  $G_{\restriction X}$ and $G'_{\restriction X}$  are isomorphic. We have:
\begin{lemma} \label{nonn-1} For every integer  $v$, $v\geq 4$, 
there are two graphs $G$ and $G'$, on the same set of $v$ vertices, which are $k$-hypomorphic for  all $k\in \{v-1,v\}$  but $G'\neq G$ and $G'\neq {\overline G}$.
\end{lemma}
\Proof Let $V:= \{0,\dots,v-1\}$, $\mathcal E:= \{\{i, i+1\}: 0\leq i<v-1\}\cup \{\{0,v-1\}\}$, $\mathcal E':= (\mathcal E \setminus \{\{0,v-1\}, \{1,2\}\})Ê\cup \{\{1,v-1\}, \{0,2\}\}$.  Let $G:= (V, \mathcal E)$ and  $G':= (V, \mathcal E')$. These  graphs are  two $v$-element cycles, $G'$ being obtained from $G$ by exchanging $0$ and $1$. Trivially, they satisfy the conclusion of the lemma. 
\endproof

With this lemma, the proof of  the first part of item $2$ is complete.

The fact that $(v, k)\in \mathcal S$ implies $k\leq \vartheta(v)$ for infinitely many $v$ is an immediate consequence of the following proposition. 

\begin{proposition} \label{nonn-3} For every integer  $v:= m+r$, where $m$ is a product of  prime powers $p$,   $p \equiv 1  \ (mod \ 4)$, and $r\in \{2,3,4\}$ 
there are two graphs $G$ and $G'$, on the same set of $v$ vertices, which are $k$-hypomorphic up to complementation for all $k$, $\vartheta (v)+1\leq k\leq v$ but $G'\neq G$ and $G'\neq {\overline G}$.
\end{proposition}

\Proof Our construction relies on the following claim.  

\begin{claim}\label{claimprime} For each integer $m$, where $m$ is a product of  prime powers $p$,   $p \equiv 1  \ (mod \ 4)$,  there is a graph $P$ on $m$ vertices such that  $P$ and $P_{-x}$ are self complementary for every $x\in V(P)$.
\end{claim}

Let $P$ be a graph satisfying the conclusion of   Claim \ref{claimprime}.  

{\bf Case 1.} $r:= 4$. In this case $\vartheta(v):= m$. Let $V$ be made of $V(P)$ and  four new elements added, say $1,2,3,4$.  Let $G$ and $G'$  be the  graphs with vertex set $V$ which coincide with $P$ on $V(P)$, the other edges of $G$ being $\{1,2\}$, $\{2,3\}$, $\{3,4\}$, $\{2,x\}$, $\{3,x\}$ for all $x\in V(P)$, the other edges of $G'$ being $\{1,3\}$, $\{2,3\}$, $\{2,4\}$, $\{2,x\}$, $\{3,x\}$ for all $x\in V(P)$. Clearly,  $G'\neq G$ and $G'\neq {\overline G}$.  Let $X\subseteq V$ with $\vert X\vert \leq 3$ and $K:= V\setminus X$. If $X\cap \{1,2,3,4\}\in \{\{1,3\}, \{2,4\}\}$ then  $G_{\restriction K}\simeq \overline {G'}_{\restriction K}$. In all other cases $G_{\restriction K}\simeq G'_{\restriction K}$. Hence $G$ and $G'$ are $k$-hypomorphic for $\vartheta(v)+1\leq k\leq v$. 

{\bf Case 2.} $r:= 3$. In this case $\vartheta (v): =m$. Let $G_1:= G_{-1}$ and $G'_1:= G'_{-1}$ where $G$, $G'$ are the graphs constructed in Case 1. Clearly $G'\neq G$ and $G'\neq {\overline G}$.  And since $G,G'$ are $k$-hypomorphic for $m+1\leq k \leq m+4$, the graphs $G_1$ and $G'_1$ are $k$-hypomorphic for $\vartheta (v)+1\leq k\leq v$.

{\bf Case 3.} $r:= 2$. In this case $\vartheta (v): =m-1$.  Let $V$ be made of $V(P)$ and  two new elements added, say $1,2$.  Let $G$ and $G'$  be the  graphs with vertex set $V$ which coincide with $P$ on $V(P)$, the other edges of $G$ being  $(2,x)$ for all $x\in V(P)$, the other edges of $G'$ being $(1,x)$ for all $x\in V(P)$. Clearly,  $G'\neq G$ and $G'\neq {\overline G}$. Let $X\subseteq V$ with $\vert X\vert  \leq 2$ and $K:= V\setminus X$. If $X\cap \{1,2\}\not = \emptyset $ then  $G_{\restriction K}\simeq \overline {G'}_{\restriction K}$. In all other cases $G_{\restriction K}\simeq G'_{\restriction K}$. Hence, $G$ and $G'$ are $k$-hypomorphic for $\vartheta(v)+1\leq k\leq v$.  \endproof \\
Claim \ref{claimprime} is an immediate consequence of the  following result.
\begin{lemma}
  Let ${\cal G}$ be the class of finite graphs $G$ of order distinct from $2$  such that $G_{-x}$ is self-complementary for every vertex $x\in V(G)$. 
  \begin{enumerate}
\item The class $\mathcal G$ coincides with the class of self-complementary vertex-transitive graphs. 
\item The class ${\cal G}$ includes Paley graphs.
\item The class ${\cal G}$ is closed under lexicographic product.
 \end{enumerate}
 \end{lemma}
  \Proof
  
$1.$ Let $G\in \mathcal G$. Let $n:=\vert V(G)\vert$.  We may suppose $n>2$. Let $x\in V(G)$. We have $d_G(x)= e(G)-e(G_{-x})$. Since $G_{-x}$ is self-complementary, $e(G_{-x})=e(\overline G_{-x})$ and, since $e(G_{-x})+e(\overline G_{-x}) = {n-1\choose 2}$, $e(G_{-x})= \frac{1}{2}{n-1\choose 2}$. Thus $d_G(x)$ does not depend on $x$, that is $G$ is regular. Since $n>2$ we have $e(G)=  \frac{1}{n-2}\sum_{x\in V(G)}e(G_{-x})$  thus $e(G)= \frac{n(n-1)}{4}$. This added to $e(G_{-x})= \frac{(n-1)(n-2)}{4}$ yields $n(n-1)\equiv 0 \ (mod \ 4)$  and $(n-1)(n-2)\equiv 0 \ (mod \ 4)$. It follows that $n\equiv 1 \ (mod \ 4)$.  As it is well-known \cite{Ke}, regular graphs of order distinct from  $2$ are reconstructible.  Thus  $G$ is self-complementary. The proof that $G$ is reconstructible yields that  for every vertex $x$,  every isomorphism from $G_{-x}$ onto $\overline G_{-x}$ is induced by an isomorphism $\varphi$ from $G$ onto $\overline G$ which fixes $x$. Hence, for a given pair of vertices $x, x'$ there is an element $\vartheta\in Aut(G)$ such that $\vartheta (x)=x'$ if and only if there is an isomorphism $\varphi: G\rightarrow \overline G$ such that $\varphi(x)=x'$. It follows that each orbit of $Aut(G)$ is preserved under all  isomorphisms from $G$ onto $\overline G$. Thus, if $A$ is a union of orbits, $G_{\restriction A}\in \mathcal G$. Since members of $\mathcal G$ have odd order, there is just one orbit, proving that $Aut(G)$ is vertex-transitive.   

Let $P$ be a self-complementary  vertex-transitive graph. Clearly $P$ is not of order $2$. Let $x\in V(P)$. Since $P$ is self-complementary,  $P_{-x}$ is isomorphic to $\overline P_{-y}$ for some $y\in V(P)$. Since $Aut(\overline P)=Aut(P)$ and $Aut(P)$ is vertex-transitive, $\overline P_{-y}$ is isomorphic to $\overline P_{-x}$. Hence, $P\in \mathcal G$. 

$2.$  Let us recall that a  {\it Paley graph} is  a graph $P_p$ whose vertices are the elements of $GF(p)$, the field of $p$ elements, with  $p\equiv 1  \ (mod \ 4)$; the  edges are all pairs  $\{x,y\}$  such that  $x-y$ is a square in  $GF(p)$. As it is well-known \cite{VW} (page 176), the  automorphism group of $P_p$  acts transitively on the edges and   $P_p$ is isomorphic to its complement. Hence, $P_{p}\in \mathcal G$.

$3.$   Let $G,H\in \mathcal G$. Their  {\it lexicographic product } is the graph $G. H$ obtained by replacing each vertex of $H$ by a copy of $G$. Formally $V(G.H):=V(G)\times V(H)$ and $E(G.H)$ is the set of pairs $\{(u,v), (u',v')\}$ such that -either $v= v'$ and $\{u,u'\}\in E(G)$, -or $\{v,v'\}\in E(H)$.
Let $x:=(u,v)\in V(G.H)$. Select an isomorphism  $\varphi $ from $G$ onto $\overline G$ which fixes $u$ and an isomorphism $\vartheta$ from $H$ onto $\overline H$ which fixes $v$. Set $\psi(u',v'):= (\varphi(u'), \vartheta (v'))$. This defines an isomorphism from $G.H$ onto $\overline {G.H}$  which fixes $x$.  Hence $G.H\in \mathcal G$
 \endproof

\begin{questions}
Does ${\cal G}$ include a graph of order $n$ whenever $n \equiv 1 \ (mod
 \  4)$ ? In particular, is there such a graph of order $21$?
\end{questions}
By Theorem \ref{k=0[4]} we have:
\begin{remark} \label{cor-fin}  Let $G$ and $G'$ be two graphs on the same set $V$ of $v$ vertices.
Let  $k$ be an integer, $1\leq k\leq v-2$, $k\equiv 0  \ (mod \ 4)$. If
$G$ and $G'$ are $(v-1)$-hypomorphic and $e(G_{\restriction K})$ and $e(G'_{\restriction K})$
have the same parity up to complementation 
  for  all $k$-element subsets $K$ of $V$, then either 
$G= G'$ or $G \in \mathcal G$.
\end{remark}
\section{Conclusion}
Let $\mathcal R$ be the set of pairs $(v,k)$ such that  two graphs on the same set are isomorphic up to complementation whenever these two graphs are $k$-hypomorphic up to complementation.

Behind Ille's problem was the question of a description of $\mathcal R$. 

This seems to be a very difficult problem.  Except the trivial inclusion $\mathcal S\subseteq \mathcal R$, the fact that some pairs like  $(5, 4)$,  $(v,v-3)$ for $v\geq 7$ requires some effort  \cite {DLPS2}. 

We prefer to point out  the following problem. 
\begin{problem}
Let $v>2$. Is $(v,k)\in \mathcal SÊ\Longleftrightarrow 4\leq k\leq \vartheta (v)$?
\end{problem}

{ \centerline  {Acknowledgements}}

We  thank J.A. Bondy,  E. Salhi and S. Thomass\'e for their helpful
comments.

\end{document}